\documentclass{amsart}
\title{Ergodicity of certain cocycles over certain interval exchanges}
\theoremstyle{plain}
\newtheorem{theorem}{Theorem}[section]
\newtheorem{lemma}[theorem]{Lemma}
\newtheorem{corollary}[theorem]{Corollary}

\newtheorem*{corollary*}{Corollary}

\DeclareMathOperator{\var}{Var}
\theoremstyle{definition}
\newtheorem{definition}[theorem]{Definition}

\newcommand{\onetorus}{\mathbb{S}^1}
\def\Z{\mathbb{Z}}
\usepackage{tikz}

\author{David Ralston}
\address{Ben Gurion University, Department of Mathematics\\ POB 653 \\ Beer Sheva, 84105\\ ISRAEL}
\email{ralston.david.s@gmail.com}

\author{Serge Troubetzkoy}
\address{Aix-Marseille University, CPT, IML, Frumam, Luminy, Case 907, F-13288 Marseille, Cedex 09, France}
\email{troubetz@iml.univ-mrs.fr}

\thanks{The first author is supported by the Center for Advanced Studies at Ben Gurion University of the Negev as well as the Israel Council for Higher Education, and was partially supported by the Erwin Schr\"{o}dinger International Institute for Mathematical Physics during preparation of this manuscript. This research is partially supported by the ANR project Perturbations.}
\date{\today}

\begin{document}
\begin{abstract}
We show that for odd-valued piecewise-constant skew products over a certain two parameter family of interval exchanges, the skew product is ergodic for a full-measure choice of parameters.
\end{abstract}
\maketitle
\section{Introduction and background}

$\Z$-valued (or more generally $G$-valued where $G$ is a locally compact group) skew products are a natural construction of infinite-measure preserving transformations using ergodic sums over a finite-measure preserving transformation.  For a thorough overview of constructing skew products over irrational rotations, see \cite{conze}.  The natural generalization of an irrational rotation is an \textit{interval exchange transformation}; recent work in studying generic skew products over generic interval exchanges may be found in \cite{chaika-hubert}, where the authors establish ergodicity of  skew products for step functions over generic interval exchanges.  We present here an alternate `hands-on' approach to prove generic ergodicity for one specific construction.

Let $X = \onetorus \times \{0,1,\ldots,k-1\}$, endowed with Lebesgue measure $\mu$ (scaled so $\mu(X)=k$), and assume that $k=1 \bmod 2$.  Let $T$ be a map on $X$ defined by
\begin{equation}
T(x,\ell) = \big((x+ \alpha) \bmod 1, (\ell + I(x)) \bmod k\big),
\label{e1}\end{equation}
where $I(x)$ is the characteristic function of an interval of length $\beta$, and $\alpha$ is irrational; $\{X, T\}$ is a $\mathbb{Z}/k\Z$-valued skew product (in fact a cyclic extension) of the irrational rotation by $\alpha$.  Let $f$ be an integer-valued function on $X$.  The skew products we will consider are given by
\[T_f(x,\ell,m) = \big((x+\alpha) \bmod 1, (\ell + I(x)) \bmod k, m+ f(x,\ell)\big).\]
Denote by $S_m(x,\ell)$ the $\mathbb{Z}$-coordinate of $T^m_f(x,\ell,0)$:
\[S_m(x,\ell) = \sum_{i=0}^{m-1}f(T^i(x,\ell)).\]
Note that $\{X \times \mathbb{Z}, T_f\}$ will \textit{not} in general itself be a skew product over rotation by $\alpha$, as $f(x,\ell)$ is not independent of $\ell$.  We assume that $f$ is of mean zero, and assume further that $f$ is piecewise constant on finitely many intervals; let $\var(f)$ be the sum over $\ell$ of the (finite) variations of $f$ restricted to each $\onetorus \times \{\ell\}$.  Purely for convenience we furthermore assume that $I$ and $f$ are right-continuous; they are defined using intervals closed on the left and open on the right.

An integer $E$ is an \textit{essential value} of our skew product if for every $A \subset X$ of positive measure, there is some $i$ such that
\[ \mu \left( A \cap T^i A \cap \{(x,\ell):S_i(x,\ell)=E\}\right)>0.\]
If $E$ is an essential value, the skew product is ergodic if and only if the skew product given by $f$ into $\mathbb{Z}/(E\mathbb{Z})$ is ergodic.

We will use \textit{Koksma's inequality}: let $P$ be a partition of $\onetorus$ into $q$ intervals of equal length, let $f$ be real-valued, of bounded variation on $\onetorus$, and suppose that $x_1$ through $x_n$ are chosen such that each interval of $P$ contains exactly one $x_m$.  Then
\[ \left| \sum_{m=1}^n f(x_m) - n \int_{\onetorus} f(x)dx \right| \leq \var(f).\]

Our interval exchanges are characterized by two choices: $\alpha$ and $\beta$.  
\begin{theorem}
Let $f$ take only odd values, and assume that not every value of $f$ is a multiple of the same number.  Then the set of $\alpha,\beta$ for which the skew product is ergodic is of full measure. 
\end{theorem} 

\section{Proof}

\begin{lemma}\label{lemmalemma}
Let $f$ take integer values (not necessarily odd) and assume that not every value of $f$ is a multiple of the same number.  Further let $\beta \in (0,1)$ be fixed, and assume there is some finite, nonzero $E \in \mathbb{Z}$ which is an essential value of the skew product $\{X \times \mathbb{Z}, T_f\}$.  Then the set of $\alpha$ for which the skew product is ergodic is of full measure.
\begin{proof}
Suppose that $\beta$ is fixed and not zero. We can construct a compact, connected translation surface $M$ and a cross-section $X$ so that the
the first return map to $X$ of the geodesic flow in the direction with slope $1/\alpha$ is $T$ given by \eqref{e1} for the parameters $\alpha,\beta$.

\begin{figure}[htb]
\center{\begin{tikzpicture}
\draw (0,0) rectangle (2,2);
\draw (2.5,0) rectangle (4.5,2);
\draw (5,0) rectangle (7,2);
\draw (0.7,-0.1) -- (0.7,0.1);
\draw (0.7,1.9) -- (0.7,2.1);
\node at (0.7, -0.3) {$\beta$};
\node at (0.35, -0.3) {$i$};
\node at (0.35, 2.3) {$iii$};

\draw (3.2,-0.1) -- (3.2,0.1);
\draw (3.2,1.9) -- (3.2,2.1);
\node at (3.2, -0.3) {$\beta$};
\node at (2.85, -0.3) {$ii$};
\node at (2.85, 2.3) {$i$};

\draw (5.7,-0.1) -- (5.7,0.1);
\draw (5.7,1.9) -- (5.7,2.1);
\node at (5.7, -0.3) {$\beta$};
\node at (5.35, -0.3) {$iii$};
\node at (5.35, 2.3) {$ii$};

\end{tikzpicture}}
\caption{The translation surface $M$ for $k=3$ and $I(s) = 1_{[0,\beta).}$ The unlabeled sides are identifies to the opposite side in the same
square, the other identifications are given by roman numbers.
The cross-section $X \times \{0,1,2\}$ consists of the bottom of the three squares.  The flow in the vertical direction corresponds to $\alpha = 0$.}
\end{figure}
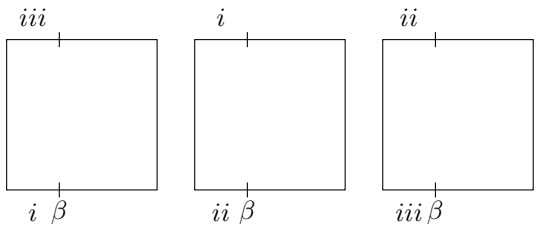

By \cite{1986}, the system $\{X, \mu, T\}$ is (uniquely) ergodic for almost every choice of $\alpha$.

Now let $X' = X \times\{0,1,\ldots,E-1\}$, with the identification
\[ (x,\ell,k) \sim (x,\ell,k+f(x,\ell) \mod E)\]
for each $(x, \ell) \in X$.  This identification corresponds to gluing together $E$ disjoint copies of $M$ via the values given by $f$, taken modulo $E$; denote this new surface by $M'$.  So long as $M'$ is connected, the results of \cite{1986} still apply, and the transformation
\[ S'(x, \ell, k) = (x+\alpha \mod 1, \ell + I(x), k+ f(x,\ell) \mod E)\]
is uniquely ergodic for almost every choice of $\alpha$.  The assumption that the values of $f$ generate $\mathbb{Z}$ exactly ensure that $M'$ is connected via B\'{e}zout's Lemma: the values taken by $f$ on each $X \times \{j\}$ do not depend on the choice of $j \in \{0,1,\ldots, E-1\}$, and there is no single common divisor for the set of values taken by $f$, so we may freely pass from one copy of $M$ to another via the values of $f$ to generate any integer value.  Ergodicity of the skew product for each $\alpha$ such that this finite system is ergodic then follows as $E$ was assumed to be an essential value of $\{X \times \mathbb{Z}, \mu \times dz, T_f\}$.
\end{proof}
\end{lemma}

The effect of Lemma \ref{lemmalemma} is to reduce our problem to the existence of a single nonzero, finite essential value for generic choice of $\beta$.  We now re-introduce the assumption that the values of $f$ are all odd (and still not multiples of the same number).  Let $\alpha$ be irrational with continued fraction expansion
\[ \alpha = [a_1,a_2,\ldots] = \cfrac{1}{a_1+\cfrac{1}{a_2+\cfrac{1}{\ddots}}}\] where each $a_m$ is a positive integer; an excellent reference for the theory of continued fractions is \cite{MR1451873}.  Denote by $p_n/q_n$ the convergents to $\alpha$, and by $\| \cdot \|$ the distance to the nearest integer.   Then it is well-known that
\begin{equation}\label{eqn - cont frac est}q_n\|q_n \alpha\| \leq \frac{1}{a_{n+1}}.\end{equation}
On $X$ we also use $\| \cdot \|$ for distance, with the convention that if $\ell \neq \ell'$, $\|(x,\ell)-(y,\ell')\|=1$.  We denote by $Q_n(T)$ the periodic approximation to $T$ given by
\[Q_n(x,\ell) = \left( x + \frac{p_n}{q_n} \bmod 1, \ell + I(x) \bmod k \right).\]

\begin{definition}
A point $x \in X$ will be called \textit{$n$-good for rational approximation} if for all $i=0,1,\ldots,kq_n-1$ we have
\[f(T^i x) = f(Q_n^i(x)), \quad I(T^i x) = I(Q_n^i x).\]  That is, as far as the functions $f$ and $I$ are concerned, through time $kq_n$ we may replace the orbit of $x$ under $T$ with the orbit of $x$ under $Q_n$.
\end{definition}

\begin{definition}
A point $x \in X$ will be called \textit{$n$-spread out} if the set $\{T^i(x)\}$, $i=0,1,\ldots,kq_n-1$, has the property that
\begin{itemize}
\item there are exactly $q_n$ points in each $\onetorus \times \{\ell\}$, and
\nopagebreak
\item for each $\ell$, there is a partition of $\onetorus \times \{\ell\}$ into disjoint intervals of length $1/q_n$ such that there is exactly one of the $T^ix$
in each partition element.
\end{itemize}
\end{definition}

\begin{lemma}
Suppose that $x$ is $n$-spread out.  Then
\[ \left|\sum_{i=0}^{kq_n-1} f(T^i x) \right|\leq \var(f).\]
\begin{proof}
The restriction of the orbit of $x$ to each $\onetorus \times \{\ell\}$ may be summed separately, and the $n$-spread out assumption allows us to use Koksma's inequality on each $\onetorus \times \{\ell\}$.
\end{proof}
\end{lemma}

Let $D=\{d_1,\ldots,d_N\}$ be the projection of all discontinuities of $f$ onto $\onetorus$ together with the discontinuities of $I(x)$.  For $n=0 \bmod 2$ define
\[A_n = \left( \onetorus \setminus \left( \bigcup_{i=0}^{kq_n-1} \bigcup_{j=1}^N \Big[d_j - k\|q_n \alpha\|- i\alpha, d_j - i \alpha \Big)\right) \right) \times \{1,2,\ldots,k\},\] while for $n=1 \bmod 2$ we use the intervals
\[ \Big(d_j - i\alpha, d_j + k\|q_n \alpha\|- i \alpha \Big].\]
\begin{lemma}\label{lemma - size of A}
Each $x \in A_n$ is $n$-good for rational approximation, and
\[\mu(A_n) \geq k\left(1 - k^2Nq_n \|q_n \alpha\|\right) \geq k \left(1- \frac{k^2N}{a_{n+1}}\right).\]
\begin{proof}
The first inequality is elementary (assume all removed intervals are disjoint), and the final inequality is simply due to \eqref{eqn - cont frac est}; the only content to prove is that $x \in A_n$ implies that $x$ is $n$-good for rational approximation.  Suppose that $n=0 \bmod 2$ so that $p_n/q_n > \alpha$.  Let $x \in A_n$; there is no $i < kq_n$ such that
\[x+i \alpha \in \left[ d_j-k\|q_n \alpha\|, d_j\right).\]
The distance between $x+ i \alpha$ and $x+ip_n/q_n$ is no larger than $k\|q_n \alpha\|$, so we cannot have
\[x+i \alpha < d_j \leq x+ i \frac{p_n}{q_n}\] for any $i,j$.  As $p_n/q_n>\alpha$, this completes the proof for $n = 0\bmod 2$.  For $n=1\bmod 2$ the process is identical, but we remove intervals from the other side of the discontinuities $d_j$, and $p_n/q_n < \alpha$.
\end{proof}
\end{lemma}
\begin{definition}
The action of $T^{kq_n}$ on $A$ is \textit{nearly-rigid} if  $\|x - T^{kq_n}(x)\| \leq k\|q_n \alpha\|$ for all $x \in A$.
\end{definition}
\begin{lemma}\label{lemma - quasirigid}
The action of $T^{kq_n}$ on $A_n$ is \textit{nearly-rigid}.
\begin{proof}
Through time $q_n$ the point $x$ orbits into the interval defining $I(x)$ some number of times.  Under $Q_n$, however, $x$ has returned exactly to the same $\onetorus$ coordinate.  Over the next $q_n$ times, the orbit of $x$ will therefore intersect this interval \textit{the same number of times} (recall that $I(x,\ell)$ is independent of $\ell$), and so on for each $q_n$ steps in the orbit.  Whatever this number of intersections is, once we have applied $Q_n$ a total of $kq_n$ times, the total number of points in these intervals must be zero modulo $k$: $Q_n^{kq_n}(x) = x$.  As $x \in A_n$, we certainly have $T^{kq_n}(x)$ belonging to the same copy of $\onetorus$ as $x$, then, and the distance in $\onetorus$ between $x$ and $T^{kq_n}(x)$ is equal to $\|kq_n \alpha\|$, which is no larger than $k\|q_n \alpha\|$.
\end{proof}
\end{lemma}
\begin{definition}
The set $A$ is \textit{nearly invariant} under $T$ if \[\mu(A \triangle T(A)) \leq 2k^2N \|q_n \alpha\|.\]
\end{definition}
\begin{lemma}
The set $A_n$ is nearly invariant under $T$.
\begin{proof}
Recall that $A_n$ is constructed by removing successive preimages of $kN$ different intervals of length $k\|q_n \alpha\|$ ($N$ such intervals in each copy of $\onetorus$).  Therefore $A_n \triangle T(A_n)$ at most consists of the first image of these intervals and the next preimage.
\end{proof}
\end{lemma}

Define \[\sigma_n(x) = \sum_{i=0}^{q_n-1} I\left( x + \frac{i}{q_n} \bmod 1\right).\]
Note that if $x \in A_n$, then
\[ \sigma_n(x) = \sum_{i=0}^{q_n-1} I(T^i x).\]
\begin{lemma}
If $x \in A_n$, $a_{n+1} \geq k$, and $\sigma_n(x)$ is relatively prime to $k$, then $x$ is $n$-spread out.
\begin{proof}
Note that $\sigma_n(x)$ is exactly the number of times through time $q_n$ that $I(Q_n^ix)=1$.  By the assumption that $x \in A_n$, this is also the number of times that $T^ix$ will orbit into this interval, and furthermore this number will be repeated for each successive length-$q_n$ segment of the orbit we consider:
\[x \in A_n \quad \Longrightarrow \quad \sigma_n(x) = \sigma_n(T^{q_n}x) = \ldots = \sigma_n(T^{(k-1)q_n}x).\]
As $\sigma_n(x)$ was assumed to be relatively prime to $k$ (i.e. $\sigma_n(x)$ generates $\mathbb{Z}/k\mathbb{Z}$), it follows that for each $i = 0,1,\ldots,q_n-1$, each of
\[\{T^{i+\ell q_n}(x)\} \quad (\ell=0,1,\ldots,k-1)\] belongs to a \textit{different} copy of $\onetorus$.  Finally, the assumption that $a_{n+1} \geq k$ implies (again via \eqref{eqn - cont frac est}) that
\[ k\|q_n \alpha\| < \frac{1}{q_n},\]
so the intervals $[x+i/q_n,x+(i+1)/q_n)$ in each circle (if $n=0 \bmod 2$; for $n=1 \bmod 2$ reverse which end is closed versus open) each contain one element of the orbit.
\end{proof}
\end{lemma}

\begin{lemma}
For all $x$, $\sigma_n(x) \in \{M, M+1\}$, where $M =[q_n \beta]$, the integer part of $q_n \beta$.
\begin{proof}
The number $M$ is the minimum number of abutting intervals of length $1/q_n$ (closed on the left, open on the right, say) which will always be completely contained within an interval of length $\beta$:
\[ \frac{M}{q_n} \leq \beta < \frac{M+1}{q_n}.\]  For any $x$, then, there are at least $M$ successive $I(x+i/q_n)=1$.  On the other hand, as $(M+1)/q_n > \beta$, no $x$ may have $\sigma_n(x) \geq M+2$.
\end{proof}
\end{lemma}

\begin{definition}
If $T^{kq_n}$ is nearly rigid and there is some $\epsilon>0$ such that $\mu(A_n) \geq \epsilon$ then $T$ is called \textit{quasi-rigid} and the $A_n$ are called \textit{quasi-rigidity sets}.
\end{definition}
\begin{corollary}\label{corollary - finite ess val exists}
Suppose that for infinitely many $n$ we have
\begin{itemize}
\item $a_{n+1} > k^2 N$,
\item $q_n = 1 \bmod 2$,
\item $\sigma_n(x)$ is relatively prime to $k$ for all $x \in X$.
\end{itemize}
Then there is a finite nonzero essential value.
\begin{proof}
The assumption that $a_{n+1} > k^2N$ implies that the $A_n$ are quasi-rigidity sets (via Lemmas \ref{lemma - size of A} and \ref{lemma - quasirigid}).  That $\sigma_n(x)$ is relatively prime to $k$ ensures that for each $x \in A_n$, $x$ is $n$-spread out, so by applying the Koksma inequality there is a uniform bound on the absolute value of the ergodic sums on $A_n$.  We therefore apply \cite[Corollary 2.6]{conze-fraczek} (utilizing that the $A_n$ are quasi-rigid and nearly invariant, which we have already established) to find an essential value (possibly zero) for the skew product; in short, as there is an upper bound on the sums from Koksma's inequality, we may pass to a sequence of subsets along which a single value is seen, and this value is therefore an essential value.   As $k q_n$ is odd and $f$ takes only odd values, it follows that for all $x \in A_n$ we must have
\[ \left|\sum_{i=0}^{kq_n-1} f(T_f^i(x)) \right| \geq 1,\] so therefore the essential value we have found in this manner is not zero.
\end{proof}
\end{corollary}

It is therefore of interest to determine when $\sigma_n(x)$ is relatively prime to $k$.  
\begin{lemma}
\label{superlacunary}
Let $\{m_i\}$ be an unbounded sequence of integers, and let $k$ be a positive integer.  Then for each residue class $j \mod k$, for almost every $\theta$ the equality
\[ [m_i \theta]=j \mod k\]
is satisfied for infinitely many $i$.
\begin{proof}
Without loss of generality, assume that $\{m_i\}$ are unbounded above, and by passing to a subsequence, we may assume that the $m_i$ are \textit{superlacunary}:
\[ \lim_{i \rightarrow \infty} \frac{m_{i+1}}{m_i} = \infty.\]  Also, without loss of generality assume $\theta \in [0,1]$, and define the random variable
\[X_i(\theta) = [m_i \theta] \mod k.\]
Suppose that $X_{i-1}(\theta)=R$, so that for some $M$ we have
\[ \theta = \frac{R+Mk}{m_{i-1}} + \frac{ \{m_{i-1} \theta\}}{m_{i-1}},\]
where $\{x\}$ denotes the fractional part of $x$.  The residue class of $[m_i \theta]$, then, is determined by the residue class of $R'$, where
\[ \theta \in \left[ \frac{R'}{m_i}, \frac{R'+1}{m_i}\right).\]
As the $\{m_i\}$ are superlacunary, the number of intervals of length $1/m_i$ within an interval of length $1/m_{i-1}$ diverges, from which it follows that
\[ \lim_{i \rightarrow \infty} \mathbb{P}\left(X_{i+1}=j | X_i \right) = \frac{1}{k}\]
for each residue class $j$.  So along this superlacunary subsequence, for generic $\theta$ the sequence $[m_i \theta]$ is uniformly distributed among the residue classes, from which the lemma trivially follows.
\end{proof}
\end{lemma}

\begin{corollary}
For almost every choice of $\alpha$, $\beta$, there are infinitely many $n$ such that such that $a_{n+1}>k^2 N$, $q_n=1 \bmod 2$, and $[q_n \beta]=1 \bmod k$.
\begin{proof}
For generic $\alpha$ there are infinitely many pairs $a_{n+1}$, $a_{n+2}$ of arbitrarily large partial quotients, and no two consecutive $q_n$, $q_{n+1}$ may be even, so the first two conditions are trivially satisfied.  The $\{q_n\}$ are an increasing sequence of integers, so by Lemma \ref{superlacunary}, for almost every $\beta$ arbitrary residue classes of $[t_m\beta]$ modulo any fixed $k$ are achieved infinitely many times.

\end{proof}
\end{corollary}
This completes the proof of ergodicity: for generic choice of $\alpha$, $\beta$ the skew product will have a nonzero essential value $E$ by Corollary \ref{corollary - finite ess val exists} (as $k$ is odd, both one and two are relatively prime to $k$).  By Lemma \ref{lemmalemma}, this suffices for generic ergodicity.

\bibliography{construction-bibfile}
\bibliographystyle{plain}

\end{document}